\documentclass[12pt, reqno]{amsart}

\usepackage{amsmath}
\usepackage{amsthm}
\usepackage{amssymb}
\usepackage{latexsym}
\usepackage{mhequ}
\usepackage{graphicx}
\usepackage{url}
\usepackage{enumerate}

\usepackage[notref,notcite]{showkeys}

\theoremstyle{plain}
\newtheorem{theorem}{Theorem}[section]
\newtheorem{lemma}[theorem]{Lemma}
\newtheorem{proposition}[theorem]{Proposition}

\newtheorem{conjecture}[theorem]{Conjecture}

\theoremstyle{definition}

\newtheorem{definition}[theorem]{Definition}

\theoremstyle{remark}
\newtheorem{remark}[theorem]{Remark}
\theoremstyle{problem}
\newtheorem{problem}[theorem]{Problem}

\numberwithin{equation}{section}

\def \prt {\partial}
\def\bone{{\bf 1}}
\def\al{\alpha}
\def\eps{\varepsilon}

\def\P{\mathbb{P}}
\def\E{{\mathcal E}}
\def\bE{{\bf E}}
\def\J{{\mathcal J}}
\def\bJ{{\bf J}}
\def\R{\mathbb{R}}

\def\BB{{\mathcal B}}
\def\MM{{\mathcal M}}

\def\wh{\widehat}

\newcommand{\comment}[1]{}

\newcommand{\eq}{\begin{equation}}
\newcommand{\en}{\end{equation}}

\def \ol{\overline}

\def\wt{\widetilde}
\def\wh{\widehat}

\begin{document}

\title{\bf Deterministic approximations of random reflectors}
\author{
{\bf Omer Angel}, {\bf Krzysztof Burdzy} and {\bf Scott
Sheffield} 
}
\address{Department of Mathematics,
University of British Columbia,
121 - 1984 Mathematics Rd,
Vancouver, BC, V6T 1Z2,
Canada}
\email{angel@math.ubc.ca}
\address{Department of Mathematics, Box 354350,
University of Washington, Seattle, WA 98195, USA}
\email{burdzy@math.washington.edu}
\address{Department of Mathematics,
Massachusetts Institute of Technology, 2-180,
77 Massachusetts Ave.,
Cambridge, MA 02139, USA }
\email{sheffield@math.mit.edu}

 \thanks{
O.~Angel's research was supported in part by NSERC and by the Sloan Foundation.
K.~Burdzy's research was supported in part by NSF Grant DMS-0906743 and
by grant N N201 397137, MNiSW, Poland. 
Scott Sheffield's research was supported in part by NSF Grant
DMS 0645585.}

\begin{abstract}
Within classical optics, one may add microscopic ``roughness'' to a macroscopically flat mirror so that parallel rays of a given angle are reflected at different outgoing angles.  Taking the limit (as the roughness becomes increasingly microscopic) one obtains a flat surface that reflects {\em randomly}, i.e., the transition from incoming to outgoing ray is described by a probability kernel (whose form depends on the nature of the microscopic roughness).

We consider two-dimensional optics (a.k.a.\ billiards) and show that {\em every} random reflector on a line that satisfies a
necessary measure-preservation condition (well established in the theory of billiards) can be approximated by deterministic reflectors in this way.
\end{abstract}

\maketitle

\baselineskip = 1.3\baselineskip

\section{Introduction }

This article addresses the question of which random reflectors
can be approximated by piecewise smooth surfaces that reflect
light according to the classical rule of specular reflection which says that the angle of reflection is equal to the
angle of incidence. One of the basic results on billiards says
that a certain measure on the space of pairs consisting  of location and angle
of reflection is preserved by every reflecting surface (see \cite[Thm. 3.1]{T} or \cite[Lemma 2.35]{CM}).
Our main result, Theorem \ref{m16.2},
shows that, except for this universal restriction that applies
to all reflecting surfaces, one can approximate in a weak sense every random
reflector by a sequence of specularly (deterministically)
reflecting surfaces.

A special case of this theorem is that one can approximate
every {\em deterministic} reflector that preserves the appropriate measure, including, for example, the reflector that reverses
the direction of each incoming ray.  To see how counterintuitive this is, imagine a higher dimensional analog of this reflector: a mirror that appears entirely the color of the observer's eyeball, because the only light rays traveling from the mirror to the eyeball are those that bounced off of the eyeball
before reaching the mirror.
This striking effect is implemented in practice
with limited accuracy (for theoretical
and practical reasons) in ``retroreflectors'' (see \cite{retro})
and reflective paint (see \cite{scotch}).

There are several sources of inspiration for our project. The
article \cite{BBCH} studies reflected Brownian motion with
inert drift. In a follow up project, Z.-Q.~Chen and the second
author plan to study the limiting situation when the diffusion
coefficient of reflected Brownian motion goes to zero. It
appears that in the limit, the particle will move along
straight lines with a random angle of reflection.

The physics literature on random reflections is quite rich. We
will not review it here; an excellent review can be found, for
example, in \cite{CPSV1}.

Some early mathematical articles that considered random
reflections were \cite{LR,E}, and more recent ones include
\cite{CP,CPSV1,CPSV1er,CPSV2,CPSV4,Fer1,Fer2,LN1,LN2}. Many of these articles are concerned
with the so called Knudsen reflection law. Knudsen proposed the
cosine reflection law, in which the direction of the reflected molecule has a
rotation-invariant distribution around the surface
normal, with density proportional to the cosine of the angle with the normal, independent of the incidence direction. Knudsen's law
is a model for gas dynamics. The same law is known as Lambert's
cosine law in optics (see \cite{CH}, pp.~147-148 or \cite{KSK}, Chap.~6).

On the technical side, our results seem to be related, at least
at the intuitive level, to the ``digital sundial'' theorem
proved by Falconer (\cite{F}, Thm.~6.9). Roughly speaking, Falconer's
theorem says that there exists a set with prescribed
projections in almost all directions.

Although our article is close to the literature on billiards at
the technical level, we will use the language of optics because
our model is much closer to this circle of ideas at the
intuitive level.

We will state our results in a rigorous way in Section \ref{main}. Section \ref{sec:det} contains short proofs of the most elementary results. The proof of our main result, Theorem \ref{m16.2}, is a multistage construction presented in Sections \ref{trans}-\ref{sec:cantor}.

\section{ The main result, conjectures and open problems }\label{main}

For a very detailed and careful presentation of the billiards
model in the plane see Chap.~2 of \cite{CM}. We will be
concerned with mirrors (walls of billiard tables) of very
special shape. They are supposed to model macroscopically flat but rough
reflecting surfaces.
The paper of Feres \cite{Fer1} contains a rigorous mathematical presentation of this physical phenomenon and detailed analysis of its fundamental properties. Our setup is slightly different from that in \cite{Fer1}.

Consider the following assumptions about a planar set $M$. These conditions contain, among other things, Assumptions A1-A4 from \cite[Sections 2.1 and 2.4]{CM}.
\medskip

(M1) $M \subset \{(x_1,x_2)\in \R^2: x_2< 0\}$.

(M2) For every $k< \infty$, the set $ \{(x_1,x_2)\in M: -k \leq
x_1 \leq k\}$ is the union of a finite number of compact $C^3$ curves $\Gamma_j$.

(M3) The curves $\Gamma_j$ intersect only at their endpoints.
Each curve either is a line segment or has non-vanishing curvature of one sign (it has no inflection points). The curves do not form cusps at the intersection points, that is, the curves form an angle different from $0$ at the point where they meet (if there is such a point).

\medskip

We will say that $M\in \MM_1$ if $M$ satisfies (M1)-(M3).

Let $D= \R \times (-\pi, 0) $ and define a $\sigma$-finite
measure $\Lambda$ on $D$ by $\Lambda(dx, d\al) = -dx \sin \al
\,d \al$.

A light ray can be represented as $(x, \al) =(x(t), \al(t))$,
where $x(t)$ is the location of the light ray at time $t$ and
$\al(t)\in(-\pi, \pi]$ is the angle between the direction of
the light ray and the positive horizontal half-axis, measured
in the counterclockwise direction from the half-axis.
Time will play no role in our arguments so it will be suppressed in the notation most of the time.
We will always assume that light rays reflect from surfaces comprising $M\in \MM_1$
according to the rule of specular reflection, that is, the
angle of incidence is equal to the angle of reflection, for every reflection.

Let $L_* :=\{(x_1,x_2) \in \R^2: x_2 =0\}$. It will be convenient to identify $L_*$ with $\R$, for example, we will consider $\Lambda$ to be a measure on $L_* \times (-\pi, 0)$.
Consider the following natural condition.

\medskip

(M4)
Suppose that $M\in \MM_1$ and
for $\Lambda$-almost all $(x,\al)\in D$, a light ray starting from $(x,\al)$ and reflecting from surfaces comprising $M$ will return to $L_*$ after a finite number of reflections.

\medskip

Condition (M4) is far from trivial; see, for example, Sec.~2.4
in \cite{CM} on accumulations of collisions. Moreover,
some light rays reflecting from some mirror
sets $M\in \MM_1$ will not return to $L_*$.
We will show that (M4) holds for a large class of sets $M$.

\medskip

(M5)
Suppose that $M\in \MM_1$. Let $\{A_k'\}_{k\geq 1}$ be the family of all connected components of the open set
$\{(x_1,x_2)\in \R^2\setminus M: x_2< 0\}$. There exists a subfamily $\{A_k\}_{k\geq 1}$ of $\{A_k'\}_{k\geq 1}$ (that is, every set $A_k$ is equal to some set $A_j'$), such that every set $A_k$ is bounded, $L_* \subset \bigcup_{k\geq 1}\prt A_k $, and the set $\{\prt A_j \cap \prt A_k \cap L_*, j,k\geq 1, j\ne k\}$ has no accumulation points in $L_*$.

\medskip

We will say that $M\in \MM_2 $ if $M\in \MM_1$ and it satisfies (M4).

\begin{proposition}\label{m13.1}
If $M\in \MM_1$ satisfies (M5) then it satisfies (M4) and, therefore, $M\in \MM_2$.
\end{proposition}

Consider some $M\in \MM_2$.
Suppose that a light ray starts from $(x_0, \al_0)$
with $x_0 \in L_*$ and $\al_0 \in (-\pi,0)$ at time 0, reflects from surfaces of $M$ and returns to $L_*$
at a time $t$, i.e., $(x(t-), \al(t-)) =(x_1, \beta')$, $x_1 \in L_*$, and $t>0$ is the smallest time with this property.
Let $\beta = \beta'-\pi$. This defines a
mapping $K : D \to D$, given by $K(x,\al) = (y,\beta)$.
Clearly, $K$ depends on $M$.

We will write $\P(x, \al; dy, d\beta)$ to denote a Markov
transition kernel on $D$, that is, for fixed $(x,\al) \in D$,
$\P(x, \al; dy, d\beta)$ is a probability measure on $D$. We
assume that $\P$ satisfies the usual measurability conditions
in all variables.

We will use $\delta_x(y)$ to denote Dirac's delta function.
Recall the transformation $K$ and let $\P_{K}$ be defined by
$\P_{K}(x, \al; dy, d\beta) = \delta_{K(x,\al)}(y,\beta)dyd\beta$. In other
words, $\P_{K}$ represents a deterministic Markov kernel, with
the atom at $K(x,\al)$.

If $\mu_n$, $n\geq 1$, and $\mu_\infty$ are non-negative
$\sigma$-finite measures on some measurable space $\Gamma$ then
we will say that $\mu_n$ converge weakly to $\mu_\infty$ if
there exists a sequence of sets $\Gamma_j$, $j\geq 1$, such
that $\bigcup_{j\geq 1} \Gamma_j = \Gamma$, $\mu_n(\Gamma_j) <
\infty$, $\mu_\infty(\Gamma_j)<\infty$ for all $n$ and $j$, and
for every fixed $j$, the sequence $\mu_n(\Gamma_j)$ converges
weakly to $\mu_\infty(\Gamma_j)$.

\begin{theorem}\label{j7.1}
(i) Consider the transformation $K: D \to D$ corresponding to any $M\in
\MM_2$. The transformation $K$ preserves measure $\Lambda$, that
is, for any $A \subset D$ with $\Lambda(A) < \infty$, we have
$\Lambda(K^{-1}(A)) = \Lambda(A)$. Moreover, $K$ is ``time
reversible'' in the sense that if $K(A_1) = A_2$ then
$K(A_2) = A_1$.

(ii) Suppose that for some sequence of sets $M_n\in \MM_2$,
corresponding transformations $K_n$, and some Markov transition kernels $ \P(x, \al; dy, d\beta)$, we have
\begin{align}\label{m16.1}
\Lambda(dx, d\al)  \P_{K_n} (x, \al; dy, d\beta)
\to
\Lambda(dx, d\al)  \P(x, \al; dy, d\beta)
\end{align}
in the sense of weak convergence
on $D^2$ as $n\to \infty$. Then $\P$ is symmetric
with respect to $\Lambda$ in the sense that for any smooth
functions $f$ and $g$ on $D$ with compact support we have
\begin{align}\label{j12.1}
\int_{D^2} f(y,\beta) \P (x, \al; dy, d\beta)
g(x,\al) \Lambda(dx, d\al)
= \int_{D^2} g(y,\beta) \P (x, \al; dy, d\beta)
f(x,\al) \Lambda(dx, d\al).
\end{align}
In particular, $\Lambda$ is invariant in the sense that
\begin{align}\label{m12.1}
\int_{D^2} f(y,\beta) \P (x, \al; dy, d\beta) \Lambda(dx, d\al)
= \int_{D} f(x,\al) \Lambda(dx, d\al).
\end{align}

\end{theorem}

See \cite[Sect.~4]{Fer1} for a similar result stated in a slightly different setting.
The first part of the theorem says that all specular
reflections are time reversible and preserve a certain measure.
For this reason, $\Lambda$ is known as the invariant measure
for the collision map in the theory of billiards (\cite{CM},
Sec.~2.12). This is related to Lambert's cosine law in optics
(see \cite{CH}, pp.~147-148 or \cite{KSK}, Chap.~6), also known as Knudsen's cosine
reflection law in the context of gas dynamics (see
\cite{CPSV1}). The second part shows that this condition can be interpreted as symmetry for a Markov kernel (see \eqref{j12.1}). This symmetry is
preserved under weak limits of Markov kernels. The next theorem, which is our main result, says that
the symmetry of the Markov kernel expressed in \eqref{j12.1} is the only condition on
a Markov kernel $\P$ necessary for the existence of
deterministic approximations of random reflections represented by $\P$.

\medskip

Recall that $\delta_x(y)$ denotes Dirac's delta function.
Suppose that the probability kernel $\P$ in Theorem \ref{j7.1} (ii)
satisfies $\P(x,\al; dy, d\beta) = \delta_x(y)dy \wt \P(x,\al; d\beta)$ for some $\wt \P$.
Heuristically, this means that the light ray released at $x$ is instantaneously reflected from a mirror located infinitesimally close to $L_*$. 
 Then \eqref{j12.1} and \eqref{m12.1} imply that for all smooth bounded functions $f$ and $g$ on $(-\pi,0)$, and almost all $x$,
\begin{align}\label{m12.2}
\int_{(-\pi,0)^2} f(\beta) \wt \P (x, \al; d\beta)
g(\al) \sin \al\, d\al
= \int_{(-\pi,0)^2} g(\beta) \wt\P (x, \al; d\beta)
f(\al) \sin\al\, d\al,
\end{align}
and
\begin{align}\label{m12.3}
\int_{(-\pi,0)^2} f(\beta) \wt \P (x, \al; d\beta) \sin \al\, d\al
= \int_{(-\pi,0)} f(\al) \sin\al\, d\al.
\end{align}

\begin{theorem}\label{m16.2}

Suppose that $\P(x,\al; dy, d\beta) = \delta_x(y)dy \wt \P(x,\al; d\beta)$ where $\wt \P$ satisfies \eqref{m12.2}.
Then there exists a sequence of sets $M_n\in \MM_2$ and corresponding
transformations $K_n$ such that
\begin{align*}
\Lambda(dx, d\al)  \P_{K_n} (x, \al; dy, d\beta)
\to
\Lambda(dx, d\al)  \P(x, \al; dy, d\beta)
\end{align*}
weakly on $D^2$ as $n\to \infty$. Moreover, $M_n$ can be chosen
in such a way that

(a) $M_n \subset \{(x_1,x_2) : -1/n < x_2 < 0\}$, and

(b) for every $\eps>0$ there exists $D_\eps \subset D$ with $\Lambda(D \setminus D_\eps) < \eps$ such that all rays starting in $D_\eps$ reflect from
$M_n$ exactly twice before returning to $L_*$.

\end{theorem}

\begin{remark}
The sets $M_n$ that we construct in the proof of Theorem
\ref{m16.2} are not connected. We believe that
the following conjecture can be proved using an
iteration of the construction used in the proof of Theorem
\ref{m16.2}.
\end{remark}

\begin{conjecture}
One can construct sets $M_n$ so that they satisfy Theorem
\ref{m16.2} except for (b), they are connected and every
light ray reflects from $M_n$ four times, except for a set of rays of $\Lambda$ measure $1/n$.
\end{conjecture}

The following problem is inspired by the ``digital sundial''
theorem of Falconer (see \cite{F}, Thm.~6.9).

\begin{problem}
Is it possible to construct sets $M_n$ so that they satisfy
Theorem \ref{m16.2} and every light ray
reflects from $M_n$ only once, except for a set of rays of $\Lambda$ measure $1/n$? 

\end{problem}

We believe that an analogue of Theorem \ref{m16.2} holds in higher dimensions.

\section{Deterministic reflections}\label{sec:det}

\begin{proof}[Proof of Proposition \ref{m13.1}.]
Assumptions A1-A4 from \cite[Sections 2.1 and 2.4]{CM} are built into (M1)-(M3). Hence, there are no accumulation points for reflections from $M$, by the results in \cite[Sect. 2.4]{CM}.

The remaining part of the proof is based on \cite[Thm. 3.1 and Sect. 7.1]{T}. We will only outline the main steps. Fix some $k$ and consider $A_k$ to be a bounded billiard table. We define the ``billiard ball map'' $T$ as follows.
Let $v(x)$ be a continuous unit tangent vector field on $\prt A_k$.
For a light ray starting from a point $x\in \prt A_k$, let $\al$ be the angle between its direction and $v(x)$. The light ray will hit $\prt A_k$ at a point $y\in \prt A_k$. Let $\beta$ be the angle formed by the direction of the light ray just after reflection at $y$ and $v(y)$. Then we let $T(x,\al)= (y,\beta)$. Let $dx$ represent the arc length measure on $\prt A_k$. By
\cite[Thm. 3.1]{T}, $T$ preserves the measure $dx\sin \al d\al$ on $\prt A_k \times (0,\pi)$.
Poincar\'e's Recurrence Theorem (see \cite[Thm. 7.4]{T})
shows that $(dx\,d\al)$-almost all rays starting at $(\prt A_k \cap L_*)\times(0,\pi)$ return to $(\prt A_k \cap L_*) \times(0,\pi)$ after a finite number of reflections. See \cite[page 116]{T}
for more details.
\end{proof}

\begin{proof}[Proof of Theorem \ref{j7.1}.]

(i) The claim that $\Lambda(K^{-1}(A)) = \Lambda(A)$
is a special case of a well known theorem, see \cite[Thm. 3.1]{T} or \cite[Lemma 2.35]{CM}. The fact that the light reflection process is time reversible implies that if $K(A_1) = A_2$ then
$K(A_2) = A_1$.

(ii)
Suppose that $M_n\in \MM_2$ and $K_n$ is the
corresponding transformation.
Consider some sets $A_1,A_2\in D$ and let $f(x,\al) = \bone_{A_1}(x,\al)$ and $g(x,\al) = \bone_{A_2}(x,\al)$.
Then, using part (i),
\begin{align*}
&\int_{D^2} f(y,\beta) \P_{K_n} (x, \al; dy, d\beta)
g(x,\al) \Lambda(dx, d\al)\\
&= \int_{D^2} \bone_{A_1}(y,\beta) \delta_{K_n(x, \al)}( dy, d\beta)
\bone_{A_2}(x,\al) \Lambda(dx, d\al)\\
&= \int_{D^2} \bone_{K_n^{-1}(A_1)}(x, \al)
\bone_{A_2}(x,\al) \Lambda(dx, d\al)\\
&= \Lambda(K_n^{-1}(A_1) \cap A_2) = \Lambda(K_n^{-1}(A_2) \cap A_1).
\end{align*}
For the same reason,
\begin{align*}
\int_{D^2} g(y,\beta) \P_{K_n} (x, \al; dy, d\beta)
f(x,\al) \Lambda(dx, d\al)
= \Lambda(K_n^{-1}(A_2) \cap A_1),
\end{align*}
so
\begin{align*}
\int_{D^2} f(y,\beta) \P_{K_n} (x, \al; dy, d\beta)
g(x,\al) \Lambda(dx, d\al)
=\int_{D^2} g(y,\beta) \P_{K_n} (x, \al; dy, d\beta)
f(x,\al) \Lambda(dx, d\al).
\end{align*}
A standard argument based on finite linear combinations of step functions and bounded convergence shows that
\eqref{j12.1} holds for smooth $f$ and $g$ with compact support and $\P$ of the form $\P_{K_n}$. This and weak convergence imply that \eqref{j12.1} also holds for Markov kernels $\P$ that satisfy \eqref{m16.1}. We obtain \eqref{m12.1} from \eqref{j12.1} by the monotone convergence theorem applied to a sequence $g_n\uparrow 1$ as $n\to \infty$.
\end{proof}

\section{Transposition reflector }\label{trans}

This section is devoted to the construction of a set of mirrors that transpose thin bundles of light of appropriate angular width. The first challenge is to guide bundles of light rays so that only an arbitrarily small amount of light is scattered in an unaccounted for way. The second challenge (easier to present in a rigorous way than the first one) is to place various sets of mirrors so that they do not interfere with each other.

Let $\ol A$ denote the closure of a set $A$.

\begin{definition}\label{def:ssf}
We will call $K:D \to D$ a {\it simple symmetric function} if there exists a countable family of rectangles $Q_k= (x_1^k, x_2^k) \times (\al_1^k, \al_2^k)\subset D$, $k\geq 1$, such that
\begin{enumerate}[(i)]
\item
$Q_k \cap Q_j = \emptyset$ for $k\ne j$, 
\item $\bigcup_{k\geq 1} \ol Q_k = D$, 
\item for any $j,k$,
either $(x_1^k, x_2^k) = (x_1^j, x_2^j)$ or $(x_1^k, x_2^k) \cap (x_1^j, x_2^j) =\emptyset$, 
\item for every $a< \infty$, there is only a finite number of $k$ such that $(x_1^k, x_2^k) \cap (-a,a) \ne\emptyset$,
\item for every $k$ there exists $\al$ such that $K(Q_k) = (x_1^k, x_2^k)\times \{\al\}$,
\item for every $k$
there exists $j$ such that
$$\int_{(\al_1^k, \al_2^k)} \sin\al d\al = \int_{(\al_1^j, \al_2^j)} \sin\al d\al,$$
$(x_1^k, x_2^k) = (x_1^j, x_2^j)$,  $K(Q_k) \subset Q_j$ and $K(Q_j) \subset Q_k$.
\end{enumerate}
We will call $\sup_k (x_2^k - x_1^k) \lor \sup_k (\al_2^k-\al_1^k)$ the {\it mesh} of $K$.
\end{definition}

Recall from Section \ref{main} that for $K:D \to D$, the kernel $\P_{K}$ is defined by
$\P_{K}(x, \al; dy, d\beta) = \delta_{K(x,\al)}(y,\beta)dyd\beta$.
Suppose that $\wt \P$ satisfies \eqref{m12.2}. Then standard arguments (see the proof of Theorem \ref{m16.2} in Section \ref{sec:cantor}) show that there exists a sequence $\{\wh K_n\}_{n\geq 1}$ of simple symmetric functions such that $\P_{\wh K_n}(x, \al; dy, d\beta) \to
\delta_x(y)dy \wt \P(x, \al; d\beta)$ weakly as $n\to \infty$.
Note that we do not claim that $\wh K_n$ arise as functions associated to reflecting sets in $ \MM_2$.
To prove Theorem \ref{m16.2}, it will suffice to show that for any simple symmetric function $\wh K_n: D\to D$ there exists a sequence of sets $M_n \in \MM_2$, $n\geq 1$, satisfying conditions (a) and (b) of
Theorem \ref{m16.2} and such that if the functions $K_n: D\to D$ correspond to $M_n$'s then $K_n - \wh K_n\to 0$ pointwise. The rest of the paper is devoted to the construction of sets $M_n$.

\bigskip

For $\rho >0$, let $\Lambda_\rho$ be the measure $\Lambda$ restricted
to $D_\rho := (-\rho, \rho) \times (-\pi, 0)$ and
note that the total mass of $\Lambda_\rho$ is $4 \rho$.

\begin{lemma}\label{j12.2}

For any $\eps_0>0 $ there exist $n_0<\infty$ and $\rho_0>0$ such that for any $n>n_0$, $\rho\in(0,\rho_0)$ and any simple symmetric function $K$ with mesh less than or equal to $1/n$ and such that $(-\rho, \rho)$ is one of the intervals $(x_1^k, x_2^k)$, the following holds.

(i) There exist a bounded set $N \in \MM_1$ and a set $D_{\rho}^*\subset D_{\rho}$ such that
$\Lambda_{\rho}(D_{\rho}^*)< 4\rho \eps_0$.

(ii) The function $K_{N}(x,\al)$ corresponding to $N$ is defined
on $D_{\rho} \setminus D_{\rho}^*$.

(iii) For all $(x,\al)\in D_{\rho} \setminus D_{\rho}^*$, we have $|K_{N}(x,\al) - K(x,\al)| \leq \eps_0$.

\end{lemma}

The function $K_N(x,\al)$ is not necessarily defined for
all $(x,\al) \in D_{\rho}$ because some light rays starting
from $L_*$ and reflected in $N$ may never come back to
$L_*$.

\begin{proof}[Proof of Lemma \ref{j12.2}]

We will omit some details of our arguments and estimates that are totally elementary but are tedious to write down.
However, we will provide now a solid justification for the approximate formulas that will form the core of the proof. Suppose that two (parts of) ellipses are fixed and serve as mirrors for light rays. Suppose that a light ray starts from $x\in L_*$ at an angle $\al \in (-\pi, 0)$, then makes two reflections from the elliptic mirrors, and then returns to $L_*$ at $y\in L_*$ and angle $\beta$. An important (but elementary) observation is that the function $(x,\al) \to (y, \beta)$ is analytic. The reason is that the equations that determine the points and angles of reflection of the light ray are quadratic with ``parameters'' that are analytic functions of $x$ and $\al$. A similar remark applies to other quantities that are functions of $x$ and $\al$, such as the location and angle of incidence when the light ray hits one of the ellipses, or the distance from the light ray to a fixed point.

\medskip
\noindent{\it Step 1}.
In this step, we will construct a pair of elliptic mirrors that interchange two very thin bundles of light rays. This step is devoted only to the construction (or, in other words, definition) of the
two mirrors. The next step will contain the proof that the mirrors reflect
the light rays in the desired directions.

We will use letters $A,B$ and $C$ (with subscripts and superscripts) to denote points in the plane. The notation $\ol{AB}$ will refer to a line segment with given endpoints;
$\ol{ABC}$ will denote an arc of an ellipse passing through the three points.

Consider some (small) $c_1>0$,
$\alpha,\beta \in (-\pi + c_1, -c_1)$ and (small) $\Delta \alpha,
\Delta \beta \in (0,1)$, such that
$\Delta \al\sin \al= \Delta\beta \sin \beta$.
Suppose that $\rho >0$ and (see Fig.~\ref{fig1}),
\begin{align*}
A_1 = (- \rho, 0),\  A_2 = (0,0),\  A_3 = (\rho,0).
\end{align*}

First assume that $\al=\beta$.
In this case, let $\E_\al = \E_\beta$ be the arc of a circle
represented in complex notation as $\E_\al=\{r e^{i\gamma}: \al-\Delta\al \leq \gamma \leq \al+\Delta\al\}$ for some $r>0$. Note that the circle is centered at $A_2$. We will consider $\E_\al$ to be a mirror. Light rays starting from $(x,\gamma)$ in $D_\al^1 := \ol {A_1 A_3} \times (\al-\Delta\al,\al+\Delta\al)$ will mostly return to the same set if $\rho$ is small. More precisely,
\begin{align*}
&
\forall c_1>0 \ \forall \eta_1>0
\ \exists \rho_0<\infty
\ \forall \rho\in(0, \rho_0)
\ \forall \al\in(-\pi+c_1,-c_1): \\
&\qquad\frac
{\Lambda_\rho(\{(x,\gamma)\in D_\al^1: K_{\E_\al}(x,\gamma)\in D_\al^1\}) }
{\Lambda_\rho(D_\al^1)} > 1-\eta_1.
\end{align*}
The case when $\al = \beta$ is rather easy so we will focus on the case $\al \ne \beta$ in the rest of the proof.

Recall notation from Definition \ref{def:ssf} and let $\wh \al^k = (\al^k_1 + \al^k_2)/2$. Consider $j$ and $k$ such that  $K(Q_k) \subset Q_j$ and $K(Q_j) \subset Q_k$. If $(\al^k_2 - \al^j_1 ) \lor (\al^j_2 - \al^k_1 ) \leq \eps_0/2$ then let $K'(x,\al) = (x, \wh \al^k)$ for $(x,\al) \in Q_k$ and $K'(x,\al) = (x, \wh \al^j)$ for $(x,\al) \in Q_j$. Otherwise, we let $K'(x,\al) = K(x,\al)$. Note that $K'$ is a simple symmetric function and $|K'(x,\al) - K(x,\al)| \leq \eps_0/2$. Hence, it will suffice to find a set $N$ and the corresponding function $K_N$ such that $|K_{N}(x,\al) - K'(x,\al)| \leq \eps_0/2$.

Consider the case when $|\alpha - \beta| \geq \eps_0/2$ and $\al-\Delta\al,\al+\Delta\al,\beta-\Delta\beta,\beta+\Delta\beta \in (-\pi,0)$.
For any angle $\gamma \in (-\pi, 0)$, let $L_\gamma$ be
the line passing through $A_2$, with the slope $\tan(- \gamma)$.
Unless stated otherwise, we will consider only points and sets below the line $L_*$.

\begin{figure}
\centering
 \includegraphics[width=15cm]{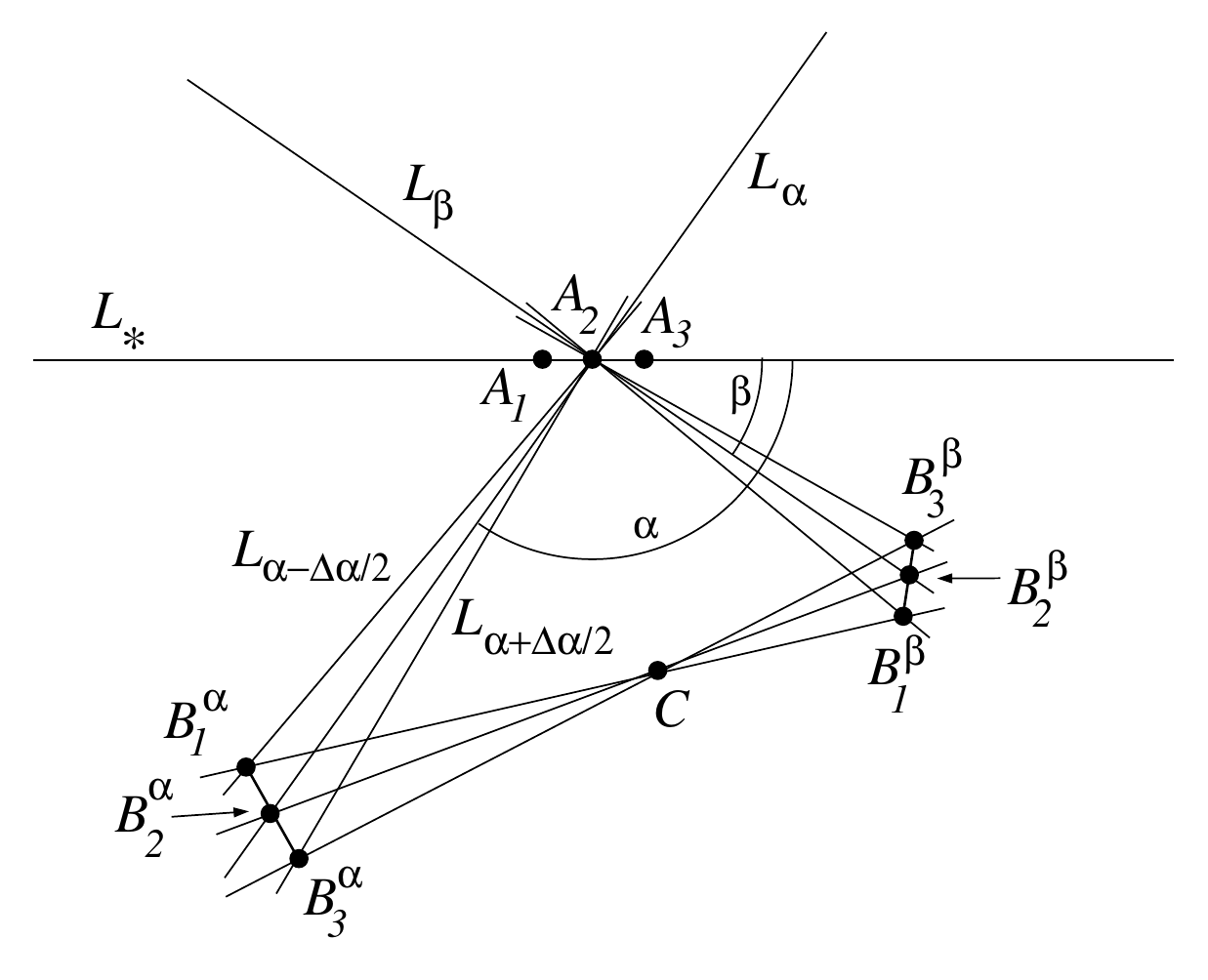}
\caption{A single light ray interchange mirror system.}\label{fig1}
\end{figure}

Let $B^\al_2 $ be a point on $L_\alpha$, and $r_\alpha := |A_2
-B^\al_2 |$. Similarly, let $B^\beta_2 $ be a point on
$L_\beta$, and $r_\beta := |A_2 -B^\beta_2 |$
(values of $r_\al$ and $r_\beta$ will be specified later). Let $L'$ be the
line passing through $B^\al_2$ and $B^\beta_2$, and let $C$ be
the point on $L'$ such that $s_1 := |C - B^\al_2| = |C -
B^\beta_2|$.

Recall that for any ellipse, if a light ray leaves one of the foci and meets a point on that ellipse, it will reflect off the ellipse and pass through the other focus.

Let $B^\al_1 \in L_{\alpha-\Delta \alpha/2}$ and $B^\al_3 \in
L_{\alpha+\Delta \alpha/2}$ be points such that $\E_\alpha :=
\ol{B^\al_1 B^\al_2 B^\al_3 }$ is an arc of an ellipse with the
foci $A_2$ and $C$ (see Fig.~\ref{fig1}). Similarly, let
$B^\beta_1 \in L_{\beta-\Delta \beta/2}$ and $B^\beta_3 \in
L_{\beta+\Delta \beta/2}$ be points such that $\E_\beta :=
\ol{B^\beta_1 B^\beta_2 B^\beta_3 }$ is an arc of an ellipse
with the foci $A_2$ and $C$.

The rest of this step is devoted to informal estimates aimed at finding values of $r_\al$ and $r_\beta$ that will define mirrors with desirable properties.
Consider the bundle $\BB$ of light rays starting at $A_2$ along lines
between $L_{\alpha-\Delta \alpha/2}$ and $L_{\alpha+\Delta
\alpha/2}$. We will find parameters of our construction such that
for the resulting mirrors $\E_\al$ and $\E_\beta$,
 ``most'' light rays in $\BB$ reflect from $\E_\alpha$, then
reflect from $\E_\beta$ and then hit $A_2$. Let $\rho_1$ be the width of $\BB$ close to
$\E_\alpha$, just before hitting $\E_\alpha$. Since $\Delta
\al$ is assumed to be small, the rays in $\BB$ are almost parallel and $\E_\al$ is almost flat. Hence,
the width of $\BB$ close to $\E_\alpha$, just after hitting
$\E_\alpha$ is $\rho_1(1 + o(\Delta\al))$. We have $\rho_1 = r_\al (\Delta
\al + o(\Delta \al))$ for small $\Delta \al$.

Since $\E_\al$ is an arc of an ellipse with foci $A_2$ and $C$, all light rays in $\BB$ will pass through $C$ and then they will hit $\E_\beta$, assuming that $\E_\beta$ is large enough.
Since $\E_\beta$ is an arc of an ellipse with foci $A_2$ and $C$, the light rays in $\BB$ reflected from $\E_\beta$ will hit $L_*$ at $A_2$.

Let $\rho_2$ be the width of $\BB$
close to $\E_\beta$, just before hitting $\E_\beta$. The width
of $\BB$ close to $\E_\beta$, just after hitting $\E_\beta$ is
$\rho_2(1+ o(\Delta\al))$. We want the bundle of light rays $\BB$ to form a cone with angle $\Delta \beta$ at the point $A_2$ on the way out. So we would like to have $\rho_2 = r_\beta (\Delta \beta +
o(\Delta \beta))$ for small $\Delta \beta$.
The point $C$ is half way between $B^\al_2$ and $B^\beta_2$ so
we would like to have $\rho_1 = \rho_2$. Therefore, we choose $r_\al$ and $r_\beta$ so that they satisfy
\begin{align*}
r_\al (\Delta \al + o(\Delta \al))
= \rho_1 = \rho_2 =
r_\beta (\Delta \beta + o(\Delta \beta)).
\end{align*}
More precisely, we choose $r_\al$ and $r_\beta$ so that
\begin{align}
\frac{r_\al}{r_\beta}
= \frac{\Delta \beta}{\Delta \al}
=\frac{\sin\al}{\sin\beta}. \label{j5.1}
\end{align}
The above formula incorporates our previously made assumption that $\Delta \al\sin \al= \Delta\beta \sin \beta$.

\medskip
\noindent{\it Step 2}. We will now argue that a pair of mirrors defined
in the previous step guides a thin bundle of light rays to the appropriate exit location and exit angle, with arbitrarily little ``waste'', for appropriate values of parameters of the construction.
Our argument is based on detailed analysis of small changes in the initial
conditions (the starting location and angle of the light ray) on the ``output,''
that is, the location and angle for the light ray exiting the lower half-plane.

Let $\bE_\al$ be the ellipse of which $\E_\al$ is a part and let $ \bE_\beta$ have the analogous meaning.
The informal estimates of Step 1 have the following rigorous version. For any $c_1,\eps>0$ there exists $\eps_1 \in(0,\eps)$ such that for all
$\alpha,\beta \in (-\pi + c_1, -c_1)$ with $|\alpha - \beta| \geq \eps_0/2$, $\Delta \alpha,
\Delta \beta \in (0,\eps_1)$, and $r_\al$ and $r_\beta$ satisfying \eqref{j5.1}, we have the following. If a light ray starts from $A_2$, follows $L_{\al + \al_1}$ with $\al_1 \in [-(1-\eps_1)\Delta\al/2, (1-\eps_1)\Delta\al/2]$ and reflects from $\E_\al$
then it will intersect $\bE_\beta$ at a point $x$ such that
\begin{align}
|x- B_2^\beta| < (1-\eps)(|B_1^\beta - B_2^\beta| \lor |B_3^\beta - B_2^\beta|). \label{s21.2}
\end{align}
Similarly,
if a light ray starts from $A_2$, follows $L_{\beta + \beta_1}$ with $\beta_1 \in [-(1-\eps_1)\Delta\beta/2, (1-\eps_1)\Delta\beta/2]$ and reflects from $\E_\beta$
then it will intersect $\bE_\al$  at a point $y$ such that
\begin{align}
|y- B_2^\al| < (1-\eps)(|B_1^\al - B_2^\al| \lor |B_3^\al - B_2^\al|). \label{s21.3}
\end{align}
It follows from \eqref{s21.2}-\eqref{s21.3} that all light rays starting from $A_2$ at angles in  $[-(1-\eps_1)\Delta\al/2, (1-\eps_1)\Delta\al/2]$
or $[-(1-\eps_1)\Delta\beta/2, (1-\eps_1)\Delta\beta/2]$
reflect from both mirrors $\E_\al$ and $\E_\beta$ and then reach $L_*$ again at $A_2$.

In the following part of the argument,
we will consider $\al, \beta, r_\al $ and $r_\beta$ to be ``fixed'' parameters satisfying \eqref{j5.1} and we will choose $\Delta \al, \Delta \beta$ and $\rho$ sufficiently small, so that certain conditions are satisfied.

We will define a number of arcs, lines and points which will be used in the next part of the proof; see Fig.~\ref{fig2}.
Consider any point $A_4 \in \ol{A_1 A_3}$ and light ray $R_1$ emanating from $A_4$ at an angle $\eta$.
Let $B_4\in \bE_\al$ and $B_6\in\bE_\beta$ denote the points where this light ray reflects in $\bE_\al$ and $\bE_\beta$. Let $R_4$ be the line between $B_4$ and $B_6$. Let $R_6$ be the line that contains this light ray after reflection in $\bE_\beta$ and let $B_8$ be the point where $R_6$ intersects $L_*$.
Let $R_2$ be the line passing through $A_2$ and parallel to $R_1$ and let $B_5$ be the intersection point of $R_2$ with $ \bE_\al$.
Let $R_3$ be the line passing through $A_2$ and $B_4$.
Let $R_5$ be the line that represents light ray $R_3$ after reflecting from $\bE_\al$.
Let $B_7$ be the point where $R_5$ intersects $ \bE_\beta$.
Let $R_7$ be the line passing through $A_2$ and $B_6$.

\begin{figure}
\centering
 \includegraphics[width=15cm]{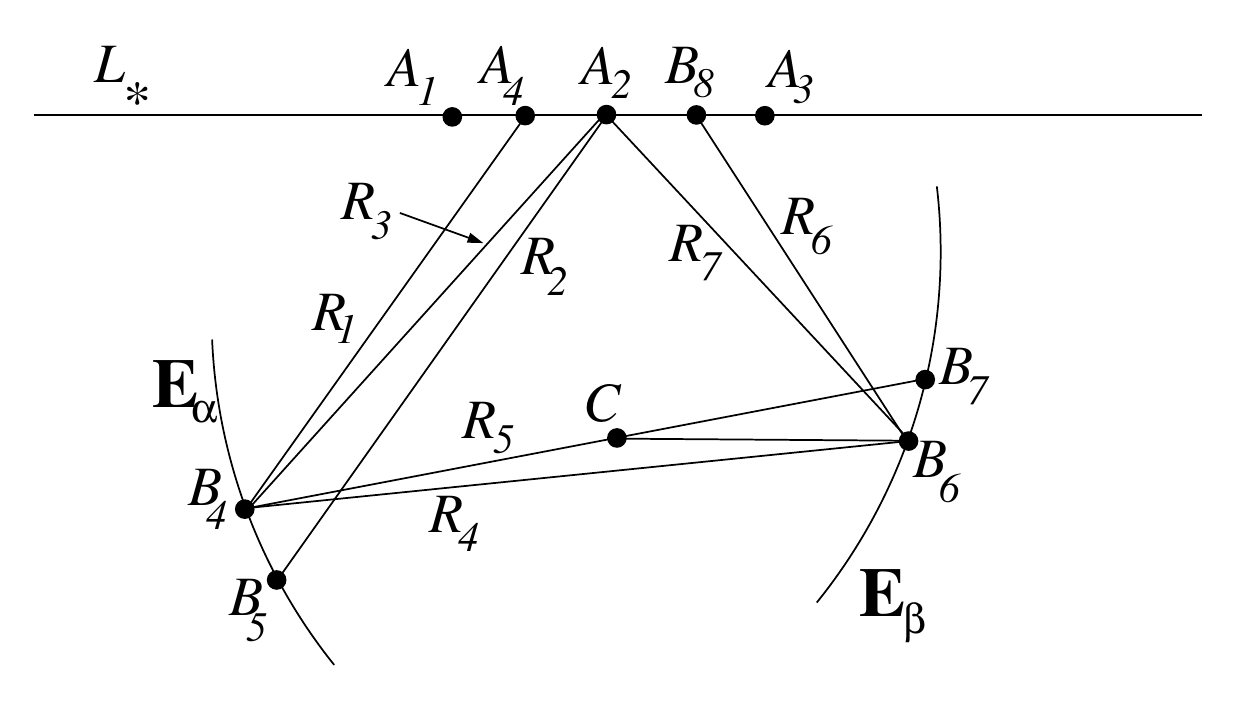}
\caption{Trajectories of various light rays.}\label{fig2}
\end{figure}

Let $\al_1$ be the angle between $\ol{A_2 B_2^\al}$ and the tangent to $\E_\al$ at $B_2^\al$. The angle $\al_1$ is a function of $\al,\beta, r_\al$ and $r_\beta$ only.
Assume that
$\eta \in  [\al-(1-\eps_1)\Delta\al/2,\al+ (1-\eps_1)\Delta\al/2]$. Then
\begin{align*}
|B_5 - B_2^\al| &\leq (1-\eps_1 + o(\eps_1))
(|B_1^\al - B_2^\al| \lor |B_3^\al - B_2^\al|)\\
& \leq
(1-\eps_1 + o(\eps_1))
\frac{r_\al (\Delta\al /2)(1+o(\Delta\al))}{\sin \al_1}.
\end{align*}
We also have
\begin{align*}
|B_4 - B_5| \leq (1 + o(\rho))
\frac{\rho \sin \al}{\sin \al_1} ,
\end{align*}
and, therefore,
\begin{align*}
|B_4 - B_2^\al| \leq
(1-\eps_1 + o(\eps_1))
\frac{r_\al (\Delta\al /2)(1+o(\Delta\al))}{\sin \al_1}
+ (1 + o(\rho))
\frac{\rho \sin \al}{\sin \al_1}.
\end{align*}
It follows that for any $c_1,\eps>0$ there exist $\eps_1 \in(0,\eps)$, $\eps_2 >0$ such that for all
$\alpha,\beta \in (-\pi + c_1, -c_1)$ with $|\alpha - \beta| \geq \eps_0/2$ and $\Delta \alpha,
\Delta \beta \in (0,\eps_1)$, $r_\al$ and $r_\beta$ satisfying \eqref{j5.1}, all 
$$\rho \in(0,\eps_2(\Delta\al \land \Delta\beta)) \text{  and  } \eta \in  [\al-(1-\eps_1)\Delta\al/2,\al+ (1-\eps_1)\Delta\al/2],$$ 
we have
$B_4 \in \E_\al$. By symmetry,
if  $\eta \in  [\beta-(1-\eps_1)\Delta\beta/2,\beta+ (1-\eps_1)\Delta\beta/2]$ then $R_1$ intersects $\bE_\beta$ at a point in $\E_\beta$.

Suppose that $\eta \in  [\al-(1-\eps_1)\Delta\al/2,\al+ (1-\eps_1)\Delta\al/2]$.
Let $\beta_1$ be the angle between $\ol{A_2 B_2^\beta}$ and the tangent to $\E_\beta$ at $B_2^\beta$. Note that $\beta_1$ is also the angle between
between $\ol{C B_2^\beta}$ and the tangent to $\E_\beta$ at $B_2^\beta$. Let $\gamma_1$ be the angle between $R_1$ and $R_3$. The angle between $R_4$ and $R_5$ is also $\gamma_1$. We have
\begin{align*}
\gamma_1 =\sin \al\, |A_4 - A_2| (1 + o(\rho) + o(\Delta \al))/r_\al \leq \rho \sin \al(1 + o(\rho) + o(\Delta \al))/r_\al.
\end{align*}
Thus
\begin{align}
|B_6 - B_7|
&= \frac{2 s_1  \sin \gamma_1 }{\sin \beta_1}  (1 + o(\rho) + o(\Delta\beta)) \label{s21.1}\\
&\leq \frac{2 s_1 \rho \sin \al (1 + o(\rho) + o(\Delta \al))}{r_\al \sin \beta_1}  (1 + o(\rho) + o(\Delta\beta))\nonumber\\
&= \frac{2 s_1 \rho \sin \al }{r_\al \sin \beta_1}  (1 + o(\rho) + o(\Delta\beta)). \nonumber
\end{align}
It follows from \eqref{s21.2}, \eqref{s21.3} and \eqref{s21.1} that for any $c_1,\eps>0$ there exist $\eps_1\in(0,\eps)$, $\eps_2 >0$ such that for all
$\alpha,\beta \in (-\pi + c_1, -c_1)$ with $|\alpha - \beta| \geq \eps_0/2$ and $\Delta \alpha,
\Delta \beta \in (0,\eps_1)$, $r_\al$ and $r_\beta$ satisfying \eqref{j5.1}, and all $\rho \in(0,\eps_2(\Delta\al \land \Delta\beta))$ we have the following. If a light ray starts from a point in $\ol{A_1 A_3}$ at an angle $\eta \in  [\al-(1-\eps_1)\Delta\al/2,\al+ (1-\eps_1)\Delta\al/2]$ then it reflects from $\E_\al$ and
then it intersects $\E_\beta$. Similarly,
if a light ray starts from a point in $\ol{A_1 A_3}$ at an angle $\eta \in  [\beta-(1-\eps_1)\Delta\beta/2,\beta+ (1-\eps_1)\Delta\beta/2]$ then it reflects from $\E_\beta$
and then it intersects $\E_\al$.

Let $a$ denote the distance between $C$ and $R_4$. We have
$$a = s_1 \sin \gamma_1 (1 + o(\rho)) = s_1 \sin \al\, |A_4 - A_2|(1 + o(\rho) + o(\Delta \al))/r_\al.$$
Note that if a light ray moved along $\ol{C B_6}$ and reflected from $\bE_\beta$ then 
it stays on $R_7$.
The angle between $R_6$ and $R_7$ is equal to the angle between $\ol{CB_6}$ and $R_4$ and, therefore, it is equal to $\gamma_2 := (a/ s_1) (1+ o(\rho)+ o(\Delta \al))$.
The following estimate for the distance between $A_2$ and $B_8$ uses \eqref{j5.1},
\begin{align*}
|B_8 - A_2| &= \frac{\gamma_2 r_\beta }{\sin \beta}
(1+ o(\Delta \al))
=
\frac{a r_\beta }{s_1\sin \beta}
(1+ o(\rho)+ o(\Delta \al))\\
&=
\frac{s_1 \sin \al\, |A_4 - A_2| r_\beta }{r_\al s_1\sin \beta}
(1+ o(\rho)+ o(\Delta \al))\\
&= |A_4 - A_2|
(1+ o(\rho)+ o(\Delta \al)).
\end{align*}
We combine this estimate with the conclusions obtained so far to see that
for any $c_1,\eps>0$ there exist $\eps_1\in(0,\eps)$, $\eps_2 >0$ such that for all
$\alpha,\beta \in (-\pi + c_1, -c_1)$ with $|\alpha - \beta| \geq \eps_0/2$ and $\Delta \alpha,
\Delta \beta \in (0,\eps_1)$, $r_\al$ and $r_\beta$ satisfying \eqref{j5.1}, all $\rho \in(0,\eps_2(\Delta\al \land \Delta\beta))$, and all
$A_4 \in \ol{A_1A_3}$ satisfying $|A_4 - A_2| \leq (1-\eps_2)\rho$
we have the following. If a light ray starts from $A_4$ at an angle $\eta \in  [\al-(1-\eps_1)\Delta\al/2,\al+ (1-\eps_1)\Delta\al/2]$ then it reflects from $\E_\al$, next it reflects from $\E_\beta$ and then it intersects $L_*$ between $A_1$ and $A_3$. Similarly,
if a light ray starts from $A_4$ at an angle $\eta \in  [\beta-(1-\eps_1)\Delta\beta/2,\beta+ (1-\eps_1)\Delta\beta/2]$ then it reflects from $\E_\beta$, next it reflects from $\E_\al$ and then it intersects $L_*$ between $A_1$ and $A_3$.

\medskip
\noindent{\it Step 3}. We will now assemble a finite family of pairs of mirrors so that they properly guide light rays entering the system at ``most'' angles and the mirrors do not interfere with one another. The lack of interference will be achieved by inductive scaling of the mirrors, that is,
making them large, so that a light ray traveling between a pair of mirrors 
takes a path far beyond all the mirrors constructed earlier in the inductive procedure.

Recall the notation from the statement of the lemma.
Let $c_1 > 0$ be so small that 
$$\Lambda_\rho \left(D_\rho \setminus \big((-\rho, \rho) \times (-\pi +c_1, -c_1)\big)\right )< 4\rho \eps_0/16 .$$
Note that $c_1$ satisfying this condition can be chosen independently of $\rho$.

Recall the simple symmetric function $K'$ defined in Step 1 and let $\{Q'_k\}_{k\geq 1} = \{(x_1^k, x_2^k) \times (\al_1^k, \al_2^k)\}_{k\geq 1}$ be the family of rectangles as in Definition
\ref{def:ssf}.
Let $\bJ$ be the family of pairs $(j,k)$ such that $j\leq k$,
$$K'((-\rho, \rho)\times (\al^j_1, \al^j_2)) \subset (-\rho, \rho)\times(\al^k_1, \al^k_2)$$ and, therefore,
$$K'((-\rho, \rho)\times (\al^k_1, \al^k_2)) \subset (-\rho, \rho)\times(\al^j_1, \al^j_2).$$ 
Note that according to Definition \ref{def:ssf} (iv), the set $\bJ$ is finite.
Let $\J $ be the set of all $j$ such that
$[\al^j_1 , \al^j_2] \cap (-\pi +c_1, -c_1) \ne \emptyset$.
For $(j,k) \in \bJ$, let 
\begin{align*}
\wh \al^j &= (\al^j_1 + \al^j_2)/2,\\
\Delta \wh \al^j &= \al^j_2 - \al^j_1,\\
\Delta \wh \al^k &= \Delta \wh \al^j \sin \wh \al^j/\sin \wh \al^k.
\end{align*}
Note that,
typically, it is not true that
$\Delta \wh \al^k = \al^k_2 - \al^k_1$ because if $(j,k) \in \bJ$ and $j\ne k$ then $(k,j) \notin \bJ$.
The quantity $\Delta \wh \al^k$ is well defined because if $(j,k) \in \bJ$ and $j\ne k$ then there is no $i$ such that $(k,i) \in \bJ$.
It is easy to see that for fixed $c_1,\eps_1>0$ we can choose $n_0$ so large that for all $n\geq n_0$, if the mesh of $K$ is smaller than $1/n$ then for all $j\in \J$ we have
$$[\wh \al^j -(1-\eps_1)\Delta \wh \al^j /2,
\wh \al^j + (1-\eps_1)\Delta \wh \al^j/2]
\subset (\al^j_1, \al^j_2).$$
Let
\begin{align*}
\Theta(\eps_1) = (-\pi +c_1, -c_1)
\cap \bigcup_{j\in \J}
[\wh \al^j -(1-\eps_1)\Delta \wh \al^j /2,
\wh \al^j + (1-\eps_1)\Delta \wh \al^j/2].
\end{align*}
For a fixed $c_1>0$,
we can make $\eps_1>0$ smaller and $n_0$ larger, if necessary, so that
\begin{align*}
\Lambda_\rho ( (-\rho, \rho) \times \Theta(\eps_1) )> 4\rho(1- \eps_0/8).
\end{align*}
We make $n_0$ larger, if necessary, so that for $n\geq n_0$ and all $j\in\J$ we have
$\Delta \wh \al^j < \eps_1$.

Choose $\eps_3 >0$ so small that
\begin{align*}
\Lambda_\rho ( (-(1-\eps_3)\rho, (1-\eps_3)\rho) \times \Theta(\eps_1) )> 4\rho(1- \eps_0/4).
\end{align*}

We will write $r(\al)$ instead of $r_\al$ and $\E(\al)$ instead of $\E_\al$, for typographical reasons.

Let $\Delta^*\al = \min_{j\in \J}\Delta \wh \al^j$.
It follows from what
we have shown in Step 2 that for $c_1,\eps>0$
there exist $n_0$, $\eps_1\in(0,\eps)$, $\eps_2 >0$ such that for all
$\rho \in(0,\eps_2\Delta^*\al)$,
$n\geq n_0$ and all $j\in\J$ we have
the following. If $(j,k) \in \bJ$ and a light ray starts from a point in $(-(1-\eps_3)\rho, (1-\eps_3)\rho)\subset L_*$ at an angle 
$$\eta \in [\wh \al^j -(1-\eps_1)\Delta \wh \al^j /2,
\wh \al^j + (1-\eps_1)\Delta \wh \al^j/2]$$
then it reflects from $\E(\wh \al^j)$, next it reflects from $\E(\wh \al^k)$ and then it intersects $L_*$ at a point  in $(-(1-\eps_3)\rho, (1-\eps_3)\rho)$.
Moreover, if a light ray starts from a point in $(-(1-\eps_3)\rho, (1-\eps_3)\rho)\subset L_*$ at an angle 
$$\eta \in [\wh \al^k -(1-\eps_1)\Delta \wh \al^k /2,
\wh \al^k + (1-\eps_1)\Delta \wh \al^k/2]$$ 
then it reflects from $\E(\wh \al^k)$, next it reflects from $\E(\wh \al^j)$ and then it intersects $L_*$ at a point  in $(-(1-\eps_3)\rho, (1-\eps_3)\rho)$.

Let $(j_1,k_1), (j_2,k_2) ,\dots$ be an arbitrary ordering of pairs in $\bJ$. Choose $r(\wh \al^{j_1})$ and $r(\wh \al^{k_1})$ so that they satisfy \eqref{j5.1}.
Suppose that $r(\wh \al^{j_m})$ and $r(\wh \al^{k_m})$ have been chosen for $m=1,2, \dots, i$. We choose
$r(\wh \al^{j_{i+1}})$ and $r(\wh \al^{k_{i+1}})$ so that they satisfy \eqref{j5.1} and they are so large that lines connecting any points in $\E(\wh \al^{j_{i+1}})$ and $\E(\wh \al^{k_{i+1}})$ do not intersect any
$\E(\wh \al^{j_{m}})$ and $\E(\wh \al^{k_{m}})$
for $m=1,2, \dots, i$.

Recall that $\bJ$ is finite and let $n_1$ be its cardinality.
For $(j,k)\in \bJ$, let $\wh D_{\rho, (j,k)}$ be the set
of light rays in $D_\rho$ such that the ray
starts from a point in $[-(1-\eps_3)\rho, (1-\eps_3)\rho]\subset L_*$ at an angle 
$$\eta \in [\wh \al^j -(1-\eps_1)\Delta \wh \al^j /2,
\wh \al^j + (1-\eps_1)\Delta \wh \al^j/2]$$ 
or at an angle 
$$\eta \in [\wh \al^k -(1-\eps_1)\Delta \wh \al^k /2,
\wh \al^k + (1-\eps_1)\Delta \wh \al^k/2].$$ 
Let $\wt D_{\rho, (j,k)}$ be the set
of light rays in $\wh D_{\rho, (j,k)}$ such that at some time the ray reflects from a set $\E(\wh \al^m)$, for some $m\ne j,k$.
If $\rho =0$, that is, if $[-(1-\eps_3)\rho, (1-\eps_3)\rho]$ is the single point $A_2$, then 
$\wt D_{\rho, (j,k)}=\emptyset $, by the claim made in the previous paragraph.
The proportion (in terms of the measure $\Lambda_\rho$) of light rays  
in $\wh D_{\rho, (j,k)}$ which belong to $\wt D_{\rho, (j,k)}$
is a continuous function of $\rho>0$, with zero limit when $\rho\downarrow 0$, that is,
$\lim_{\rho\downarrow 0} \Lambda_\rho(\wt D_{\rho, (j,k)})/ \Lambda_\rho(\wh D_{\rho, (j,k)})=0$.
It follows that for sufficiently small $\rho>0$, we have
$\Lambda_\rho(\wt D_{\rho, (j,k)}) < 4\rho\eps_0/(16 n_1)$.
Let $\wt D_\rho =\bigcup_{(j,k)\in \bJ} \wt D_{\rho, (j,k)}$. 
We see that for fixed $r(\wh \al^{j})$'s, we can make $\rho>0$ so small that $\Lambda_\rho(\wt D_\rho) < 4\rho\eps_0/16$.

Let $N = \bigcup_{j \in \J} \E(\wh \al^j)$.
Our arguments have shown that the function $K_N$ satisfies the assertions listed in the lemma, with function $K'$ in place of $K$, and $\eps_0/2$ in place of $\eps_0$. We have observed in Step 1 that this is sufficient for the proof of the lemma.
\end{proof}

\section{Light reflector cells in Cantor set holes}\label{sec:cantor}

The system of mirrors constructed in the previous section is very ``inefficient'' in that it has a large diameter compared to the length of the segment $(-\rho,\rho)$ of $L_*$ where the light enters and exits the system. In this section, we will build a ``compact'' version of the reflector using scaling and a Cantor-like construction.
The following result is very similar to Lemma \ref{j12.2} except that we use smaller mirrors.

\begin{lemma}\label{j12.3}

Fix any $\rho, \rho_1 >0$, $a\in\R$ and let $\Lambda_\rho$ be the measure $\Lambda$ restricted
to $D_\rho := (a-\rho,a+ \rho) \times (-\pi, 0)$.
The total mass of $\Lambda_\rho$ is $4 \rho$. For
any $\eps_0>0 $ there exists $n_0<\infty$ such that for any $n>n_0$ and any
simple symmetric function $K$ with parameter $n$ such that $[a-\rho,a+ \rho]$ is one of the intervals $[x_1^k, x_2^k]$, the following holds.

(i) There exist a compact set $N \in \MM_1$ and a set $D_{\rho}^*\subset D_{\rho}$ such that
$\Lambda_{\rho}(D_{\rho}^*)< 4\rho \eps_0$.

(ii) The function $K_{N}(x,\al)$ corresponding to $N$ is defined
on $D_{\rho} \setminus D_{\rho}^*$.

(iii) For all $(x,\al)\in D_{\rho} \setminus D_{\rho}^*$, we have $|K_{N}(x,\al) - K(x, \al)| \leq \eps_0$.

(iv) $N\subset (a-\rho, a+\rho) \times (-\rho_1, 0)$.

\end{lemma}

\begin{proof}

By translation invariance, it is enough to discuss the case $a=0$.

Recall that we consider $(-\rho, \rho)$ to be a subset of the horizontal axis $L_*$.
Let $N_0$ be a new name for the set satisfying Lemma \ref{j12.2} and let $M$ be the closure of the convex hull of $N_0\cup (-\rho, \rho)$.
Since $c_1$ in the proof of Lemma \ref{j12.2} is strictly positive, $M \cap L_* = [-\rho, \rho]$.
Let $r_0$ be the diameter of $M$ and let $M_r = \{x\in \R^2: \exists y \in M \text{ such that } x= (r/r_0) y\} $, i.e., $M_r$ is a dilation of $M$.

We will define sets $M_r^j$. Each of these sets will be a horizontal shift of $M_r$, that is, $M_r^j = M_r + (b,0)$ for some $b\in\R$ depending on $j$ and $r$. For any such set, we let $\wt M^j_r = M_r^j \cap L_*$. Note that the length of $\wt M^j_r$ is $2(r/r_0)\rho$.

The sets $M_r^j$ will be grouped into a countable number of families, with all sets in 
one family having the same size. We will pack smaller sets among the bigger sets
in a tight way, as much as possible. The intersections of the sets $M_r^j$
with $L_*$ will form a pattern qualitatively similar to holes in the Cantor set. See Fig.~\ref{fig7}.
\begin{figure}
\centering
 \includegraphics[width=10cm]{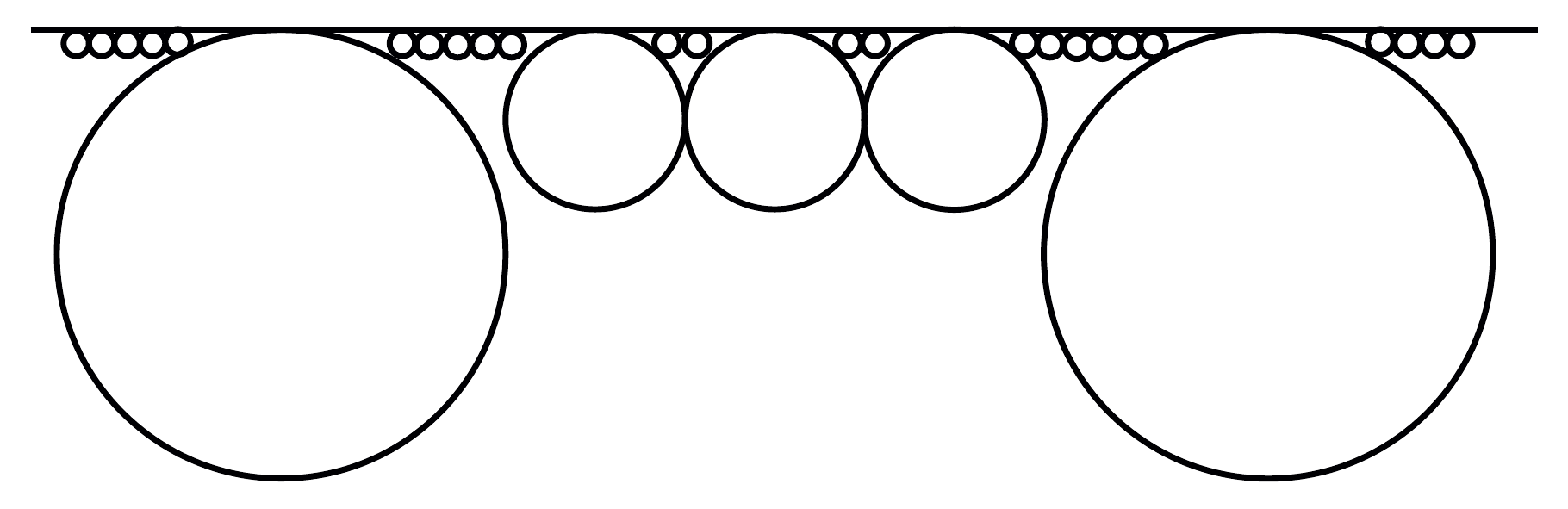}
\caption{Sets $M^j_{r_k}$ of three ``generations.'' The sets $M_r^j$
are convex by definition. They are pictured as circles
but  in fact they are not discs.}\label{fig7}
\end{figure}

For $r\in(0,1)$, let $n_r$ be the maximum number $k$ of disjoint sets $M_r^1, M_r^2, \dots, M_r^k$ such that
$\bigcup_{1\leq j \leq k} M^j_r \subset \{(x_1,x_2): |x_1| < 1\}$. It is easy to see that for some $c_2,r_*>0$ and all $r\in(0,r_*)$ we have $r n_r > c_2 $. It follows that for some $c_3>0$, for every $r\in(0,r_*)$ one can find a family of disjoint sets $M_r^1, M_r^2, \dots, M_r^k$ such that $\bigcup_{1\leq j \leq k} M^j_r \subset \{(x_1,x_2): |x_1| < 1\}$ and the total length of
$\bigcup_{1\leq j \leq k} \wt M^j_r$ is greater than $c_3$.

Choose $r_1 \in (0,\rho_1)$ and a family of disjoint sets $\{M_{r_1}^j\}_{1\leq j \leq j_1}$ such that $\bigcup_{1\leq j \leq j_1} M^j_{r_1} \subset \{(x_1,x_2): |x_1| < \rho\}$ and the total length of
$\bigcup_{1\leq j \leq j_1} \wt M^j_{r_1}$ is greater than $c_3 \rho$.

We proceed by induction. Let $\lambda_1$ be the Lebesgue measure of the set $[-\rho, \rho] \setminus
\bigcup_{1\leq j \leq j_1} \wt M^j_{r_1}$.
We find
$r_2 \in (0,r_1)$ and a family of sets $\{M_{r_2}^j\}_{1\leq j \leq j_2}$ such that 
\begin{enumerate}[(i)]
\item $\bigcup_{1\leq j \leq j_2} M^j_{r_2} \subset \{(x_1,x_2): |x_1| < \rho\}$, 
\item all sets in the family $\{M_{r_k}^j\}_{k=1,2; 1\leq j \leq j_k}$ are disjoint, and 
\item the total length of
$\bigcup_{1\leq j \leq j_2} \wt M^j_{r_2}$ is greater than $c_3 \lambda_1/2$.
\end{enumerate}

The general inductive step is the following. Suppose that we have defined $r_k>0$ and families
$\{M_{r_k}^j\}_{1\leq j \leq j_k}$ for $k=1,2, \dots, i$.
Let $\lambda_i$ be the Lebesgue measure of the set $[-\rho, \rho] \setminus \bigcup_{1\leq k \leq i}
\bigcup_{1\leq j \leq j_k} \wt M^j_{r_k}$.
We find
$r_{i+1} \in (0,r_i)$ and a family of sets $\{M_{r_{i+1}}^j\}_{1\leq j \leq j_{i+1}}$ such that 
\begin{enumerate}[(i)]
\item $\bigcup_{1\leq j \leq j_{i+1}} M^j_{r_{i+1}} \subset \{(x_1,x_2): |x_1| < \rho\}$, 
\item all sets in the family $\{M_{r_k}^j\}_{1\leq k \leq i+1; 1\leq j \leq j_k}$ are disjoint, and 
\item the total length of
$\bigcup_{1\leq j \leq j_{i+1}} \wt M^j_{r_{i+1}}$ is greater than $c_3 \lambda_i/2$.
\end{enumerate}

Let $T_r^j$ be the linear transformation that maps $M$ onto $M_r^j$ and let $N_r^j = T_r^j(N_0)$.
Let $N_* = \bigcup_{k \geq 1}
\bigcup_{1\leq j \leq j_k} N^j_{r_k}$.
It is easy to see that the Lebesgue measure of the set $ [-\rho, \rho] \setminus \bigcup_{k \geq 1}
\bigcup_{1\leq j \leq j_k} \wt M^j_{r_k}$ is zero.
It is also clear that
$N_*\subset \{(x_1,x_2): |x_1| \leq \rho\}$.

Let $\wh D \subset D_\rho$ be such that $\Lambda_\rho(\wh D) \geq 4\rho(1-\eps_0)$ and for all $(x,\al)\in \wh D$, we have $|K_{N_0}(x,\al) - K(x, \al)| \leq \eps_0$. Let $\wh D_r = \{(x,\al): \exists (y,\al) \in \wh D \text { such that } x=(r/r_0) y\}$. If $M_r^j = M_r + (b^j_r,0)$ then we let $\wh D_r^j  = \{(x,\al): \exists (y,\al) \in \wh D \text { such that } x= y + (b^j_r,0)\}$.
Let $\Lambda_r^j$ be the measure $\Lambda$ restricted
to $(-(r/r_0)\rho + b^j_r, (r/r_0)\rho+b^j_r) \times (-\pi, 0)$. By the scaling and shift
invariance of reflections, we see that for every $r>0$ and $j$ we have $\Lambda_r^j(\wh D_r^j) \geq (r/r_0)4\rho(1-\eps_0)$ and for all $(x,\al)\in \wh D_r^j$, we have $|K_{N_*}(x,\al) - K(x, \al)| \leq \eps_0$. Let $D_\rho^* = D_\rho \setminus \bigcup_{k \geq 1}
\bigcup_{1\leq j \leq j_k} \wh D^j_{r_k}$. Then $\Lambda_{\rho}(D_{\rho}^*)< 4\rho \eps_0$, and
for all $(x,\al)\in D_{\rho} \setminus D_{\rho}^*$, we have $|K_{N_*}(x,\al) - K(x, \al)| \leq \eps_0$.

The set $N_*$ satisfies all assertions in the lemma except that it does not belong to $\MM_1$ because it consists of an infinite number of curves (arcs of ellipses) and hence it does not satisfy condition (M2) of Section \ref{main}.
To address this problem, we let $N = \bigcup_{1\leq k \leq k_0} \bigcup_{1\leq j \leq j_k} N^j_{r_k}$
for some $k_0<\infty$.
It is easy to see that $N$ and the accompanying function
$K_N$ satisfy all the conditions stated in the lemma if $k_0$ is sufficiently large.
\end{proof}

\begin{proof}[Proof of Theorem \ref{m16.2}]

Recall the definition of a simple symmetric function from the beginning of Section \ref{trans}.
Suppose that $\wt \P$ satisfies \eqref{m12.2}.
Recall that $\int_{-\pi}^0 \sin\al d\al = -2$, let $m\geq 1$ be an integer and let $\gamma_k^m $ be defined by $-\int_{-\pi}^{\gamma_k^m} \sin\al d\al = k 2^{-m}$ for $k=0,1, \dots, 2^{m+1}$.
If we take $f(\al) = \bone_{(\gamma^m_n, \gamma^m_{n+1}]}(\al)$
and $g(\al) = \bone_{(\gamma^m_j, \gamma^m_{j+1}]}(\al)$ then
\eqref{m12.2} yields
\begin{align}\label{o2.1}
\int_{(\gamma^m_j, \gamma^m_{j+1}]} \int_{(\gamma^m_n, \gamma^m_{n+1}]}\wt \P(x,\al; d\beta)\sin\al d\al
=
\int_{(\gamma^m_n, \gamma^m_{n+1}]} \int_{(\gamma^m_j, \gamma^m_{j+1}]}\wt \P(x,\al; d\beta)\sin\al d\al.
\end{align}
For a fixed integer $-\infty < k < \infty$ let
\begin{align*}
a(n,j) =
-\int_{[k/m, (k+1)/m]}\int_{(\gamma^m_n, \gamma^m_{n+1}]} \int_{(\gamma^m_j, \gamma^m_{j+1}]}\wt \P(x,\al; d\beta)\sin\al d\al dx.
\end{align*}
It follows from \eqref{o2.1} that $a(n,j) = a(j,n)$.
Let $\beta^m_{n,0} = \gamma^m_n$ and let $\beta^m_{n,i}$ be defined by
\begin{align*}
-\int_{(k/m, (k+1)/m]}\int_{\gamma^m_n}^{\beta_{n,i}^m} \sin\al d\al dx =
\sum_{0\leq j\leq i-1} a(n,j)
\end{align*}
for $i=1, \dots, 2^{m+1}$. Note that $\int_{(-\pi,0]}\wt \P(x,\al; d\beta) =1$ so
\begin{align*}
\sum_{0\leq j\leq 2^{m+1}-1} a(n,j) =
-\int_{[k/m, (k+1)/m]}\int_{(\gamma^m_n, \gamma^m_{n+1}]} \sin\al d\al dx,
\end{align*}
and, therefore,
$\beta^m_{n, 2^{m+1}} = \gamma^m_{n+1}$.

Let $\wh \beta_{n,i}^m = (\beta_{n,i}^m + \beta_{n,i+1}^m)/2$ and $\wh K_m(x,\al) = (x,\wh \beta_{n,j}^m)$ for $(x,\al)\in (k/m, (k+1)/m] \times (\beta_{j,n}^m, \beta_{j,n+1}^m]$.
Since $a(n,j) = a(j,n)$, we have 
$$\int_{(\beta_{j,n}^m, \beta_{j,n+1}^m)} \sin\al d\al = \int_{(\beta_{n,j}^m, \beta_{n,j+1}^m)} \sin\al d\al.$$ 
This shows that $\wh K_m(x,\al)$ is a simple symmetric function according to Definition \ref{def:ssf}.

Recall that for $K:D \to D$, the kernel $\P_{K}$ is defined by
$\P_{K}(x, \al; dy, d\beta) = \delta_{K(x,\al)}(y,\beta)dyd\beta$.
It is elementary to check that
$\P_{\wh K_n}(x, \al; dy, d\beta) \to
\delta_x(y)dy \wt \P(x, \al; d\beta)$ weakly on $D$ as $n\to \infty$.
Note that we do not claim that $\wh K_n$ arise as functions associated to reflecting sets $M\in \MM_2$.

To prove Theorem \ref{m16.2}, it will suffice to show that there exists a sequence of sets $M_n \in \MM_2$, $n\geq 1$, satisfying conditions (a) and (b) of
Theorem \ref{m16.2} and such that if the functions $K_n: D\to D$ correspond to $M_n$'s then $K_n - \wh K_n\to 0$ pointwise.

Fix some $\wh K_n$ and let $(x_1^k, x_2^k)$ be one of the corresponding intervals as in Definition \ref{def:ssf}.
Let $D_k = (x_1^k, x_2^k) \times (-\pi,0)$.
According to Lemma \ref{j12.3} and its proof,
there exist a compact set $N_k \in \MM_1$ and a set $D_k^*\subset D_k$ such that
$\Lambda(D_k^*)< \eps_0 2^{-k}$, and the following assertions hold.
\begin{enumerate}[(i)]
\item
The function $K_{N_k}(x,\al)$ corresponding to $N_k$ is defined
on $D_k \setminus D_k^*$.
\item For all $(x,\al)\in D_k \setminus D_k^*$, we have $|K_{N_k}(x,\al) - \wh K_n(x, \al)| \leq \eps_0$, and the light ray with the starting position and direction represented by $(x,\al)$ reflects exactly twice from $N_k$ before returning to $L_*$.
\item $N_k\subset (x_1^k, x_2^k) \times (-1/n, 0)$.
\end{enumerate}

Note that for $j\ne k$, either the sets $N_j$ and $N_k$ do not intersect or they are identical. Moreover, the construction presented in Lemma \ref{j12.3} allows us to assume that light rays that start in $D_k \setminus D_k^*$ do not leave $(x_1^k, x_2^k) \times (-\infty, 0)$. In this sense, the families of mirrors $N_k$ do not interfere with one another.
Let $M^*_n = \bigcup_{k\geq 1} N_k$ and let $K^*_n(x,\al)$
be the function on $D$ corresponding to $M^*_n$. Let $D^* = \bigcup_{k\geq 1} D_k^*$ and note that $\Lambda(D^*) \leq \eps_0$. We have
\begin{enumerate}[(i)]
\item
The function $K^*_n(x,\al)$ is defined
on $D \setminus D^*$.
\item For all $(x,\al)\in D \setminus D^*$, we have $|K^*_n(x,\al) - \wh K_n(x, \al)| \leq \eps_0$, and the light ray with the starting position and direction represented by $(x,\al)$ reflects exactly twice from $M^*_n$ before returning to $L_*$.
\item $M^*_n\subset \R \times (-1/n, 0)$.
\end{enumerate}

It is easy to see that $M_n^* \in \MM_1$ for every $n$ but it is not necessarily true (actually unlikely) that $M_n^* \in \MM_2$. We can instead consider sets
\begin{align*}
M_n = M_n^* \cup \{(x_1,x_2): x_2 = -1/n\}
\cup \bigcup_k \{(x_1, x_2): x_1 = x_1^k \text{  or  } x_2^k, -1/n\leq x_2 \leq 0\}.
\end{align*}
It is easy to see that functions $K_n$ corresponding to $ M_n$ satisfy all properties (i)-(iii) listed above. Moreover, $M_n$ satisfy condition (M5) of Section \ref{main} so $ M_n \in \MM_2$, according to Proposition \ref{m13.1}.
This completes the proof of the theorem.
\end{proof}

\bigskip
\noindent{\bf Acknowledgment}. We are grateful to Persi
Diaconis and John Sylvester for the most helpful advice.
We thank the anonymous referee for many suggestions 
for improvement.


\begin{thebibliography}{bme}

\bibitem{BBCH} K.~Burdzy, R.~Bass, Z.~Chen and M.~Hairer,
    Stationary distributions for diffusions with inert drift
{\it Probab. Theory Rel. Fields \bf 146} (2010) 1--47.

\bibitem{CH}
S.~Chandrasekhar,
{\it Radiative transfer}. Dover Publications, Inc., New York, 1960.

\bibitem{CM} N.~Chernov and R.~Markarian, {\it Chaotic
    billiards}. Mathematical Surveys and Monographs, 127.
    American Mathematical Society, Providence, RI, 2006.

\bibitem{CP} F.~Comets and S.~Popov,
Ballistic regime for random walks in random environment with unbounded jumps and Knudsen billiards {\it Ann. Inst. H. Poincar\'e Probab. Statist.} (2012) (to appear).

\bibitem{CPSV1} F.~Comets, S.~Popov, G.~Sch\"utz and
    M.~Vachkovskaia, Billiards in a general domain with random reflections
    	{\it Arch. Rat. Mech. Anal. \bf 191}, (2009) 497--537.

\bibitem{CPSV1er} F.~Comets, S.~Popov, G.~Sch\"utz and
    M.~Vachkovskaia, Erratum: Billiards in a general domain with random reflections {\it Arch. Ration. Mech. Anal. \bf 193}, (2009), 737--738.

\bibitem{CPSV2} F.~Comets, S.~Popov, G.~Sch\"utz and
    M.~Vachkovskaia, Quenched invariance principle for the
    Knudsen stochastic billiard in a random tube {\it Ann. Probab. \bf 38} (2010) 1019--1061.

\bibitem{CPSV4} F.~Comets, S.~Popov, G.~Sch\"utz and
    M.~Vachkovskaia, Knudsen gas in a finite random tube: transport diffusion and first passage properties. 
{\it J. Stat. Phys. \bf 140} (2010) 948--984. 

\bibitem{E} S.~Evans, Stochastic billiards on general tables.
    {\it Ann. Appl. Probab. \bf 11}, 419--437 (2001).

\bibitem{F} K.~Falconer, {\it Fractal geometry. Mathematical
    foundations and applications. Second edition}. John Wiley
    \& Sons, Inc., Hoboken, NJ, 1990.

\bibitem{Fer1}
R.~Feres, 
Random walks derived from billiards. {\it Dynamics, ergodic theory, and geometry}, 179--222,
{\it Math. Sci. Res. Inst. Publ., \bf 54}, Cambridge Univ. Press, Cambridge, 2007. 

\bibitem{Fer2}
R.~Feres and H.-K.~Zhang, The spectrum of the billiard Laplacian of a family of random billiards. {\it J. Stat. Phys. \bf 141}, (2010) 1039--1054.

\bibitem{KSK} R.~Klette, K.~Schl\"uns and A.~Koschan, {\it
    Computer vision. Three-dimensional data from images}.
    Springer-Verlag Singapore, Singapore, 1998.

\bibitem{LR} S.~Lalley and H.~Robbins, Stochastic search in a
    convex region. {\it Probab. Theory Relat. Fields \bf 77},
    99--116 (1988).

\bibitem{LN1}
M.~Lapidus and R.~Niemeyer,
Towards the Koch snowflake fractal billiard: computer experiments and mathematical conjectures. {\it Gems in experimental mathematics}, 231--263,
{\it Contemp. Math., \bf 517}, Amer. Math. Soc., Providence, RI, 2010. 

\bibitem{LN2}
M.~Lapidus and R.~Niemeyer,
Families of Periodic Orbits of the Koch Snowflake Fractal Billiard. 
Math ArXiv 1105.0737

\bibitem{T} S.~Tabachnikov, {\it Geometry and billiards}.
    Student Mathematical Library, 30. American Mathematical
    Society, Providence, RI; Mathematics Advanced Study
    Semesters, University Park, PA, 2005.

\bibitem{retro} Wikipedia, {\it Retroreflector}, 
\url{http://en.wikipedia.org/wiki/Retroreflector}
Online; accessed 4-March-2012.

\bibitem{scotch} Wikipedia, {\it Scotchlite},
\url{http://en.wikipedia.org/wiki/Scotchlite}
Online; accessed 4-March-2012.



\end{thebibliography}
\end{document}